\documentclass[12pt,letterpaper,leqno]{article}

\pdfoutput=1

\usepackage{graphics,epsfig,graphicx,color}
\graphicspath{{../Figures/}}
\usepackage[caption = false]{subfig}
\usepackage{float}
\usepackage{capt-of}

\usepackage{amssymb,latexsym}

\usepackage{setspace}

\usepackage{amsmath}

\usepackage{mathrsfs,amsfonts}

\usepackage{amsthm,amsxtra}

\usepackage[final]{showkeys}

\usepackage{stmaryrd}

\usepackage{algorithm}
\usepackage{algorithmicx}
\usepackage{algpseudocode}

\usepackage[openbib]{currvita}

\newif\ifPDF
\ifx\pdfoutput\undefined
	\PDFfalse
\else
	\ifnum\pdfoutput > 0
		\PDFtrue
	\else
		\PDFfalse
	\fi
\fi

\ifPDF
	\usepackage{pdftricks}
	\begin{psinputs}
		\usepackage{pstricks}
		\usepackage{pstcol}
		\usepackage{pst-plot}
		\usepackage{pst-tree}
		\usepackage{pst-eps}
		\usepackage{multido}
		\usepackage{pst-node}
		\usepackage{pst-eps}
	\end{psinputs}
\else
	\usepackage{pstricks}
\fi

\ifPDF
	\usepackage[debug,pdftex,colorlinks=true, 
	linkcolor=blue, bookmarksopen=false,
	plainpages=false,pdfpagelabels]{hyperref}
\else
	\usepackage[dvips]{hyperref}
\fi

\usepackage{fancyhdr}

\newtheorem{theorem}{Theorem}[section]

\newtheorem{lemma}{Lemma}[section]
\newtheorem{remark}{Remark}[section]

\newtheorem{proposition}{Proposition}[section]

\newenvironment{keywords}
{\noindent{\bf Key words.}\small}{\par\vspace{1ex}}
\newenvironment{AMS}
{\noindent{\bf AMS subject classifications 2010.}\small}{\par}

\newcommand{\be}{\begin{equation}}
\newcommand{\ee}{\end{equation}}
\newcommand{\ba}{\begin{array}}
\newcommand{\ea}{\end{array}}

\newcommand{\bea}{\begin{eqnarray*}}
\newcommand{\eea}{\end{eqnarray*}}
\newcommand{\bean}{\begin{eqnarray}}
\newcommand{\eean}{\end{eqnarray}}

\def\tilde{\widetilde}

\def\x{{\bf x}}

\def\cydot{\leavevmode\raise.4ex\hbox{.}}

\title{A global stability estimate for the photo-acoustic inverse problem in layered media}

\author{
	Kui Ren\thanks{
		Department of Mathematics and the Institute of Computational Engineering
and Sciences (ICES), University of Texas, Austin, TX 78712, USA;
		\href{mailto:ren@math.utexas.edu}{ren@math.utexas.edu}
	}
	\and
	Faouzi Triki\thanks{
		Laboratoire Jean Kuntzmann, UMR CNRS 5224, Universit\'{e} Grenoble-
Alpes, 700 Avenue Centrale, 38401 Saint-Martin-d'H\`{e}res, France;
		\href{mailto:faouzi.triki@univ-grenoble-alpes.fr}{faouzi.triki@univ-grenoble-alpes.fr}
	}
}

\date{May 11, 2017}
\begin{document}

\maketitle



\begin{abstract}
This paper is concerned with the stability issue in determining absorption and diffusion coefficients in photoacoustic imaging. Assuming that the medium is layered and the acoustic wave speed is known we derive global H\"{o}lder stability estimates of the photo-acoustic inversion. These results show that the reconstruction is stable in the region close to the optical illumination source, and deteriorate exponentially far away. Several experimental pointed out that the resolution depth of the photo-acoustic modality is about tens of millimeters. Our stability estimates confirm these observations and give a rigorous quantification of this depth resolution.
\end{abstract}


\begin{keywords}
Inverse problems, wave equation, diffusion equation, Lipschitz  stability, multiwave inverse problems.
\end{keywords}


\begin{AMS}
	35R30, 35J15, 35L05, 92C55
\end{AMS}


\section{Introduction}
\label{SEC:Intro}

Photoacoustic imaging (PAI)~\cite{AmBrGaJu-SIAM12,AmKaKi-JDE12,Bal-IO12,KuKu-HMMI10,LiWa-PMB09,Scherzer-Book10,Wang-Book09} is a recent hybrid imaging modality that couples diffusive optical waves with ultrasound waves to achieve high-resolution imaging of optical properties of heterogeneous media such as biological tissues.

In a typical PAI experiment, a short pulse of near infra-red photons is radiated into a medium of interest. A part of the photon energy is absorbed by the medium, which leads to the heating of the medium. The heating then results in a local temperature rise. The medium expanses due to this temperature rise. When the rest of the photons leave the medium, the temperature of the medium drops accordingly, which leads to the contraction of the medium. The expansion and contraction of the medium induces pressure changes which then propagate in the form of ultrasound waves. Ultrasound transducers located on an observation surface, usually a part of the surface surrounding the object, measure the generated ultrasound waves over an interval of time $(0, T)$ with $T$ large enough. The collected information is used to reconstruct the optical absorption and scattering properties of the medium.

Assuming that the ultrasound speed in the medium is known, the inversion procedure in PAI proceeds in two steps. In the first step, we reconstruct the initial pressure field, a quantity that is proportional to the local absorbed energy inside the medium, from measured pressure data. Mathematically speaking, this is a linear inverse source problem for the acoustic wave equation~\cite{AgKuKu-PIS09,AgQu-JFA96,AmBrJuWa-LNM12,BuMaHaPa-PRE07,CoArBe-IP07,FiHaRa-SIAM07,Haltmeier-SIAM11,HaScSc-M2AS05,Hristova-IP09,KiSc-SIAM13,KuKu-EJAM08,Kunyansky-IP08,Nguyen-IPI09,PaSc-IP07,QiStUhZh-SIAM11,StUh-IP09,Tittelfitz-IP12}. In the second step, we reconstruct the optical absorption and diffusion coefficients using the result of the first inversion as available internal data~\cite{AmBoCaTaFi-SIAM08,AmBoJuKa-SIAM10,BaUh-IP10,BaUh-CPAM13,MaRe-CMS14,NaSc-SIAM14,QiSa-SIAM15,ReGaZh-SIAM13}.

In theory, photoacoustic imaging provides both contrast and resolution. The contrast in PAI is mainly due to the sensitivity of the optical absorption and scattering properties of the media in the near infra-red regime. For instance, different biological tissues absorbs NIR photons differently. The resolution in PAI comes in when the acoustic properties of the underlying medium is independent of its optical properties, and therefore the wavelength of the ultrasound generated provides good resolution (usually submillimeter).

In practice, it has been observed in various experiments that the imaging depth, i.e. the maximal depth of the medium at which structures can be resolved at expected resolution, of PAI is still fairly limited, usually on the order of millimeters. This is mainly due to the limitation on the penetration ability of diffusive NIR photons: optical signals are attenuated significantly by absorption and scattering. The same issue that is faced in optical tomography~\cite{Arridge-IP99}. Therefore, the ultrasound signal generated decays very fast in the depth direction. 

The objective of this work is to mathematically analyze the issue of imaging depth in PAI. To be more precise, assuming that the underlying medium is layered, we derive a stability estimate that shows that image reconstruction in PAI is stable in the region close to the optical illumination source, and deteriorates exponentially in the depth direction. This provides a rigorous explanation on the imaging depth issue of PAI.

In the first section we introduce the PAI model and give the main global stability estimates in Theorem~\ref{mainglobalinversion}. Section~\ref{SEC:Results} is devoted to the acoustic inversion, we derive observability inequalities corresponding to the internal data generated by well chosen laser illuminations. We also provide an observability inequality from one side for general initial states in Theorem~\ref{observability}. In section~\ref{SEC:Acous}, we solve the optical inversion and show weighted stability estimates of the recovery of the optical coefficients from the knowledge of two internal data. Finally, the main global stability estimates are obtained by combining stability estimates from the acoustic and optical inversions.

\section{The main results}
\label{SEC:Results}

\begin{figure}[!ht]
\centering
\includegraphics[width=0.5\textwidth]{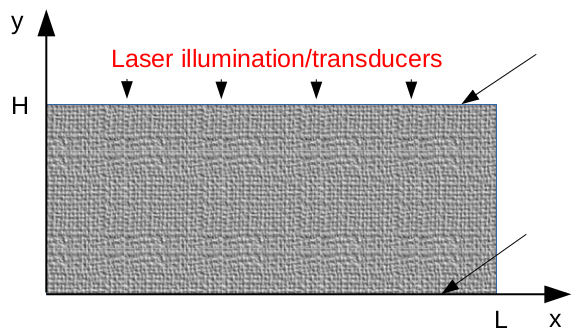}
\put(-10,40){$\Gamma_0$}
\put(-10,110){$\Gamma_m$}
\put(-150,35){$\Omega$}
\caption{The geometry of the sample.}
\label{FIG:Setup}
\end{figure}

In our model we  assume that the laser source and the ultrasound transducers are on the same side of the sample $\Gamma_m$; see Figure~\ref{FIG:Setup}. This situation is quite realistic since in applications only a part of the boundary is accessible and in the exiting prototypes a laser source acts trough a small hole  in the transducers. We also assume that the optical parameters $(D, \mu_a)$, similar to the acoustic speed $c$, only depend on the variable $y$ following the normal direction to $\Gamma_m$.  We further consider the optical parameters $(D(y), \mu_a(y))$ within
the set
\bea
\mathcal O_M =\{(D, \mu) \in C^3([0, H])^2; \;\;
D>D_0,\, \mu>\mu_0; \;\;  \|D\|_{C^3},    \|\mu\|_{C^3}\leq M\}, 
\eea
where $D_0>0, \mu_0>0$ and $M> \max(D_0, \mu_0)$ are
fixed real constants.

The propagation of the optical wave in the sample is modeled by the following diffusion equation
\bean\label{maineqO}
\left \{ \ba{llllccc} 
-\nabla \cdot D(y)\nabla u(\x) +\mu_a(y) u(\x) = 0  
&\x\in \Omega, \\
u(\x)=g(\x)&\x\in \Gamma_m, \\
 u(\x) =0 &\x\in \Gamma_0, \\
u(0,y) = {u(L, y)} & y\in (0,H), 
\ea
\right.
\eean
where $g$ is the laser illumination,
 $D$  and $\mu_a$ are respectively the  diffusion  and absorption coefficients.
The part of the boundaries $\Gamma_j$
are given by 
\bea
\Gamma_m\;=\; (0,L)\times \{y=H\}, &\Gamma_0\;=\; (0,L)\times \{y=0\},
\eea
and  $\partial_\nu$ is the derivative along $\nu$, the unit normal vector pointing 
outward of $\Omega$. 
We note that $\nu$ is everywhere defined except at the vertices of $\Omega$ and we denote by $\Gamma_p$  the complementary of $\overline{\Gamma_0}\cup\overline{ \Gamma_m}$ in
$\partial \Omega$.\\

We follow the approach taken in several papers~\cite{BaRe-IP11,BaRe-CM11,BaUh-CPAM13} and 
consider two laser illuminations $g_j, j=1,2$. Denote $u_j, j=1, 2, $ 
the corresponding  laser intensities. \\

Let
 \bea V\,:=\,\{ v \in H^1(\Omega); v(0, y) = v(L,y), \; y\in (0,H); u=0 \, \textrm{ on } \Gamma_0 \}.\eea

We further assume that $g\in V_{\Gamma_m}$, where 
 
 \bea V_{\Gamma_m}\,:=\,\{ v|_{\Gamma_m} \in H^1(\Omega); v(0, y) = v(L,y), \; y\in (0,H); 
 u=0 \, \textrm{ on } \Gamma_0 \}.\eea

Then, there exists a unique solution $u \in V$ satisfying the system \eqref{maineqO}. 
The proof  uses  techniques developed in \cite{Grisvard-Book85}. The first step is to show that the set 
\bea V_0\,:=\,\{ v \in H^1(\Omega); v(0, y) = v(L,y), \; y\in (0,H); u=0 \, \textrm{ on } 
\Gamma_0 \cup \Gamma_m\},\eea
 is a closed sub space of $H^1(\Omega)$, using a specific trace theorem 
 for regular curvilinear polygons (Theorem 1.5.2.8, page 50 in~\cite{Grisvard-Book85}).  Then, applying the 
 classical  Lax-Milligram for elliptic operators in Lipschitz domain  gives the existence 
 and uniqueness of solution to the system \eqref{maineqO}.
\begin{remark}
Using an explicit characterization of the trace theorem obtained in   \cite{Grisvard-Book85} one can derive 
the  optimal  local regularity  for $g \in H^{\frac{1}{2}}(\Gamma_m)$  that guarantees  
the existence and uniqueness of solutions to the system \eqref{maineqO} (see 
also  \cite{AmCh-DPDE17}). 
\end{remark}

For simplicity, we will further consider $g_j = \varphi_{k_j}, j=1,2$,
where $k_1<k_2$, and  $\varphi_k(x),  k\in \mathbb N, $ is the Fourier  
orthonormal basis of $L^2(0,L)$ satisfying $-\varphi_k^{\prime\prime}(y) =
\lambda_k^2 \varphi_k(y)$, with $\lambda_k =\frac{2k\pi}{L}, k\in \mathbb Z$. 
 Direct calculation gives  $\varphi_k(x) = \frac{1}{\sqrt L} e^{i\lambda_k x},  k\in \mathbb Z$.\\

We  assume that point-like ultrasound transducers, located on an
 observation surface $\Gamma_m$, are used to detect the values of the 
pressure $p(\x, t)$, where $\x \in \Gamma_m$ is a detector location and $t \geq 0$ is
 the time of the observation. We also assume that the speed of sound  in the sample
occupying $\Omega = (0,L)\times (0,H)$, is a smooth function and 
depends only on the vertical variable $y$, that is,
 $c= c(y)>0$. Then, the  following model is known 
to describe correctly the propagating pressure 
wave $p(\x, t)$ generated by the
photoacoustic effect 
\bean\label{maineq}
\left \{ \ba{llllccc} 
\partial_{tt}p(\x,t) =c^2(y) \Delta p(\x,t)  &\x\in \Omega, t\geq 0, \\
\partial_\nu p(\x,t)+\beta \partial_t p(\x,t) =0 &\x\in \Gamma_m, t\geq 0,\\
 p(\x,t) =0 &\x\in \Gamma_0, t\geq 0,\\
p((0,y),t) = p((L, y),t) & y\in (0,H), t\geq 0,\\
p(\x,0)= f_0(\x),\; \partial_{t}p(\x,0) =f_1(\x),& \x\in \Omega,
\ea
\right.
\eean
Here  $\beta>0$ is the damping coefficient, and $f_j(\x), j= 0, 1,$ are
 the initial values of the acoustic pressure, which one needs to find in order to 
determine the optical parameters of the sample. \\

\begin{remark}
Notice that in most existing works in  photoacoustic imaging 
the initial state  $f_0(\x)$ is  given by  $ \mu_a(\x) u(\x)$, and is 
assumed to be compactly supported  inside $\Omega$, while  the 
initial speed  $f_1(\x)$ is  zero everywhere \cite{KuKu-EJAM08, QiStUhZh-SIAM11, StUh-IP09}. 
The compactly support assumption on $f_0(\x)$  simplifies the analysis  of  the inverse 
source initial-to-boundary problem and is necessary for almost all the existing uniqueness
and stability  results \cite{AmChTr-PAMS16,AmChTr-arXiv16, Isakov-Book02,Yamamoto-IP95}. 
Meanwhile the assumption is clearly in contradiction with the 
 fact that  $f_0(\x)$  coincides  with  $ \mu_a(\x) u(\x)$ everywhere. 
 We will show in 
 section \ref{opticalinversion} that  $ \mu_a(\x) u(\x)$  is not only not 
 compactly supported, 
 it is also exponentially concentrated around the part of the boundary
  $\Gamma_m$ where we applied the laser illumination.  
In our model  the initial speed $\partial_{t}p(\x,0) =f_1(\x)$
can be considered as the  correction of the photoacoustic effect
generated by the heat at $\Gamma_m$. \\
\end{remark}

The following stability estimates are the main results of the paper,
obtained by combining stability estimates  from the acoustic and optical 
inversions. 

\begin{theorem} \label{mainglobalinversion}
Let $(D,\mu_a)$,\,  $(\widetilde D, \tilde \mu_a)$
in $\mathcal O_M$, and $k_i, \, i=1,2$ be two distinct integers.
Let $c(y) \in W^{1,\infty}(0, {H})$ with $ 0<c_m\leq c^{-2} (y)$
and set $\theta = \sqrt{\|c^{-2}\|_{L^\infty}}$. 
 Denote $u_{k_i}, \, i=1,2$
 and $\tilde u_{k_i}, \, i=1,2 $  the solutions  to 
 the system~\eqref{equk} for $g_{i} = \varphi_{k_i}, i=1,2$, with 
coefficients 
 $(D,\mu_a)$ and  $(\widetilde D, \tilde \mu_a)$ respectively.
 Assume that $D(H)= 
 \widetilde D(H)$,  $D^\prime(H)= 
 \widetilde D^\prime(H)$, $\mu_a^\prime(H)= 
 \widetilde \mu_a^\prime(H),\,$ $k_1<k_2$, and $k_1$ is large enough. \\

  Then,  for  $T> 2\theta H$,  there exists a constant $C>0$ 
 that only depends on  $\mu_0, D_0, k_1, k_2, M, L, $ and 
 $  H,$ such that 
 the following stability estimates  hold.
 \bea
\| \underline u_{m}^2(\mu_a- \tilde \mu_a) \|_{C^0} \leq \\ 
  C
\left( \sum_{i=1}^2 \int_{0}^{T}
 \left(
\frac{C_M}{T-2\theta H}  +\beta \right)\|\partial_{t}p_{i}- 
\partial_{t}\tilde p_{i}\|^2_{L^2(\Gamma_m)} +
 \|\partial_xp_i - \partial_x\tilde p_i \|^2_{L^2(\Gamma_m)}
 dt\right)^{\frac{1}{4}},
\eea
 and 
\bea
\| \underline u_{m}^2(D- \tilde D) \|_{C^0} \leq \\ 
C
\left( \sum_{i=1}^2\int_{0}^{T}
 \left(\frac{C_M}{T-2\theta H}  +\beta \right)\|\partial_{t}p_{i}- 
\partial_{t}\tilde p_{i}\|^2_{L^2(\Gamma_m)} +
 \|\partial_xp_i - \partial_x\tilde p_i \|^2_{L^2(\Gamma_m)}dt
\right)^{\frac{1}{4}},
\eea
where
\bea
C_M = H e^{\int_0^H  c^{2}(s)|\partial_y (c^{-2}(s))|ds } 
(c^{-2}(H) +\beta^2),\\
 \underline u_m(y) = \frac{D^{\frac{1}{2}}(H)}{D^{\frac{1}{2}}(y)}
 \frac{\sinh(\kappa_m^{\frac{1}{2}} y)}
 { \sinh(\kappa_m^{\frac{1}{2}}H)},\;  \kappa_m = \min_{0\leq y\leq H}
 \left( \frac{(D^{\frac{1}{2}})^{\prime\prime}}{D^{\frac{1}{2}}}+ \frac{\mu_a}{D}+
\lambda_k^2\right).
\eea

 \end{theorem}
Since the function $\underline u_m(y)$ is exponentially decreasing  between the value $1$
on $\Gamma_m$ to the value $0$ on $\Gamma_0$, the stability estimates in Theorem \ref{mainglobalinversion} shows that the resolution deteriorate  exponentially in the depth direction
far from $\Gamma_m$. 
 
\section{The acoustic inversion}
\label{SEC:Acous}

The data obtained by the point detectors located on the surface 
$\Gamma_m$ are represented by the function
\bea
p(\x,t)= d(\x,t) \qquad \x\in \Gamma_m, \;  t  \geq 0.
\eea
Thus, the first inversion in photoacoustic imaging  is to find, using the data 
$d(\x, t)$ measured by transducers, the initial value $f_0 (\x)$ at $t = 0$ of the 
solution $p(\x, t)$ of \eqref{maineq}.  We will also recover the initial speed $f_1(\x)$ 
inside $\Omega$, but  we will not use it in the second inversion. \\ 
 
 We first focus on the direct problem and prove existence and uniqueness 
 of the acoustic problem~\eqref{maineq}. 
 Denote {by} $L^2_c(\Omega)$ the  Sobolev  space of square 
 integrable functions with weight $\frac{1}{c^2(y)}$. Since the speed 
 $c^2$ is lower and upper bounded,  the norm corresponding to this weight is 
 equivalent to the  classical norm  of $L^2(\Omega)$.
  Let  
 \bea
 V= \{p \in H^1(\Omega); \; p(0, y) = p(L, y),\;   y\in (0,H); p = 0 \;
 \textrm{on} \; \Gamma_0 \},
 \eea
 and consider in $V\times L^2_c(\Omega) $  the unbounded linear operator $A$
 defined by 
 \bea
 A (p,q)= (q, c^2 \Delta p),\;D(A) = \{ (p,q) \in V\times V;\;  \Delta p \in L^2(\Omega); \;
 \partial_\nu  p+ \beta q = 0 \;\textrm{on} \; \Gamma_m \}.
 \eea
 
 We have the following existence and uniqueness result. 
 \begin{proposition} \label{preg}For $(f_0, f_1) \in D(A)$, 
 the problem~\eqref{maineq} has a unique solution $p(x,t)$  satisfying 
 \bea
 (p, \partial_t p)\in  C\left([0,+\infty),  D(A) \right) \cap 
 C^1\left([0,+\infty), V \times L^2_c(\Omega)\right).
  \eea
 \end{proposition}
 \proof 
 There are various methods for proving well-posedness of evolution problems:
variational methods, the Laplace transform method 
 and  the semi-group method. Here, we will  consider the semi-group method~\cite{TuWe-Book09}, and  prove
 that the operator $A$ is m-dissipative on the Hilbert space $V\times L^2_c(\Omega)$.\\

 Denote by $\langle \cdot,  \cdot\rangle$ the scalar product in $V\times L^2_c(\Omega)$,
 that is,  for $(p_i, q_i) \in V\times L^2_c(\Omega)$ with $i=1, 2$,
 \bea
 \langle (p_1, q_1), (p_2, q_2)\rangle &=& \int_\Omega \nabla p_1 \nabla  \overline 
 p_2 d\x + \int_\Omega  q_1 \overline q_2 \frac{d\x}{c^2}.
 \eea
 Now let  $ (p,q)\in D(A)$. We have 
 \bea
 \langle A (p, q), (p, q)\rangle   &=& \int_\Omega \nabla q\nabla  \overline 
 pd\x + \int_\Omega  \Delta p \overline q d\x.
 \eea
 Since $\Delta p\in L^2(\Omega)$ and $\partial_\nu  p+ \beta q = 0 \;\textrm{on} \;
  \Gamma_m$, applying Green formula leads 
 to
 \bea
 \langle A (p, q), (p, q)\rangle &=& \int_\Omega \nabla q\nabla  \overline 
 pd\x  - \int_\Omega \nabla \overline 
 q\nabla   p d\x -  \beta \int_{\Gamma_m}  |q|^2 d\sigma(\x).
 \eea
 Consequently 
 \bea
 \Re(\langle A (p, q), (p, q)\rangle) &=& -  \beta \int_{\Gamma_m}|q|^2 d\sigma(\x).
 \eea
 Therefore the operator $A$  is dissipative.  The fact that $0$ is in the resolvent
 of $A$ is straightforward. Then  $A$  is m-dissipative and hence, it is the generator of a 
 strongly continuous semigroup of contractions~\cite{TuWe-Book09}. Consequently, for 
 $(f_0, f_1) \in D(A)$ there exists  a unique strong solution to the problem~\eqref{maineq}.  
 \endproof

 Now, back to the inverse problem of reconstructing the initial data $(f_0, f_1)$.
 We further assume that the initial data is generated by a finite number  of Fourier
 modes, that is 
  \bean \label{freg}
 f_j(x,y)&=& \sum_{|k| \leq N} f_{jk}(y) \varphi_k(x) \qquad 
 (x, y)\in \Omega\qquad j=0, 1,
\eean
with $N$ being a fixed positive integer. \\
 
 As it was already remarked in many works, this linear initial-to-boundary
 inverse problem is 
 strongly related to boundary observability of the source from the set 
 $\Gamma_m$ (see for instance~\cite{LiMa-Book72,StUh-IP09, TuWe-Book09, Zuazua-CEAUCT01}). 
 We will emphasize on the links between our findings  and known
 results  in this context later. 
  Here we will use a 
 different approach taking 
 advantage of the fact that the wave speed $c(y)$ only depends in the vertical 
 variable $y$.
 
 Since $p(\x)$ is $L$-periodic in the $y$ variable, it has the following
 discrete Fourier decomposition
 \bea
 p(x, y)&=& \sum_{|k| \leq N} p_k(y,t) \varphi_k(x) \qquad 
 (x,  y)\in \Omega.
\eea

One can check  that  $p_k(y,t)  \varphi_k(x)$  is exactly   
the solution to the problem~\eqref{maineq} with initial data 
$(f_{0k}(y) \varphi_k(x), f_{1k}(y) \varphi_k(x))$.  
Precisely, if $\lambda_k = \frac{2k\pi}{L}$,
the functions $p_k(y,t) $ satisfy the
following one dimensional wave equation 
\bean\label{mainpy}
\left \{ \ba{llllccc} 
\frac{1}{c^2(y)} \partial_{tt}p(y,t) = \partial_{yy} p(y,t) - \lambda_k^2 p(y,t),  &
y\in (0,H), t\geq 0, \\
\partial_y p(H,t)+\beta \partial_t p(H,t) =0 & t\geq 0,\\
 p(0,t) =0 &t\geq 0,\\
p(y,0)= f_{0k}(y),\; \partial_{t}p(y,0) = f_{1k}(y),& y\in (0,H),
\ea
\right.
\eean
 Next, we will focus on the
boundary observability problem of the initial data $f_k$ at  the extremity 
$y = H$.  Taking advantage of the fact that the equation
is one dimensional  we will derive a boundary observability
inequality with a sharp constant.  
Define $E(t)$ the total energy  of the system~\eqref{mainpy} by
\bean\label{energ}
E(t) &=& \int_0^H  \left(c^{-2}(y)  |\partial_{t}p(y,t)|^2  +
 |\partial_{y}p(y,t)|^2 +\lambda_k^2 |p(y,t)|^2 \right) dy. 
 \eean
 Multiplying the first equation in the system~\eqref{mainpy}
  by $\partial_t p(y,t) $ and integrating
over $(0,H)$ leads to 
\bean \label{energderivative}
E^\prime(t) = - \beta  |\partial_{t}p(H,t)|^2 \qquad \textrm{for} \quad
t \geq 0.
\eean
Consequently, $E(t)$ is  a non-increasing function, and  the decay
is clearly related to the magnitude  of the dissipation on the
boundary $\Gamma_m$.

It is well know that the system~\eqref{mainpy} has 
a unique solution. Here we establish an estimate 
of the continuity constant. 
\begin{proposition} \label{pr1}
Assume that $c(y) \in W^{1,\infty}(0,{H})$ with 
$ 0<c_m\leq c^{-2} (y)$. 
Then, for any  $T>0$ we have 
\bea
\beta^2 \int_0^T  |\partial_{t}p_k(H,t)|^2 dt \leq
 ((C_m^1+C_m^2\lambda_k) T+C_m^3)  
E_k(0),
\eea
for $k\in \mathbb N$, where 
\bea
E_k(0) = \int_0^H  \left(c^{-2}(y)  |f_{1k}(y)|^2  +
 |f_{0k}^\prime(y)|^2 +\lambda_k^2 |f_{0k}(y)|^2 \right) dy,\\
C_m^1 = (1+H c^{-2}(H))^{-1}
\left( 1+(1+\frac{H}{c_m})\| c^{-2}\|_{ W^{1,\infty}}\right),\\
  C_m^3 = (1+H c^{-2}(H))^{-1}
\left(1+ 2H\| c^{-2}\|_{ L^{\infty}}^{1/2}\right),\\
C_m^2 = H(1+H c^{-2}(H))^{-1}.
\eea
\end{proposition}

\begin{proposition} \label{pr2}
Assume that $c(y) \in W^{1,\infty}(0, {H})$ with $ 0<c_m\leq c^{-2} (y)$.
Let  $\theta = \sqrt{\|c^{-2}\|_{L^\infty}}$  and $T> 2\theta H$. Then
the following inequalities hold
\bea
\lambda_k^2\int_0^H  |f_{0k}(y)|^2 dy \leq
\left(\frac{C_M}{T-2\theta H}  +
\beta \right) \int_{0}^{T}
 |\partial_{t}p_k(H,t)|^2 dt\\  + \lambda_k^2\int_{0}^{T} |p_k(H,t)|^2 dt,
\eea
for $k\in \mathbb N^*$,
\bea
\int_0^H  c^{-2}(y)  |f_{1k}(y)|^2 + |f_{0k}^\prime(y)|^2 dy   \leq 
\left(\frac{C_M}{T-2\theta H}  +\beta \right) 
\int_{0}^{T}
 |\partial_{t}p_k(H,t)|^2 dt\\  + \lambda_k^2 \int_{0}^{T} |p_k(H,t)|^2 dt,
\eea
for $k\in \mathbb N$, with 
\bea
C_M = H e^{\int_0^H  c^{2}(s)|\partial_y (c^{-2}(s))|ds } 
(c^{-2}(H) +\beta^2).
\eea
\end{proposition}
The proofs of these  results are given in the Appendix. 
The main result of this section is the following.
\begin{theorem} \label{mainacoustic}
Assume that $c(y) \in W^{1,\infty}(0,1)$ with $ 0<c_m\leq c^{-2} (y)$,
and $f_0, f_1$ have a finite Fourier expansion \eqref{freg}. 
Let  $\theta = \sqrt{\|c^{-2}\|_{L^\infty}}$  and $T> 2\theta H$. Then
\bea
\int_\Omega  |\nabla f_0(\x)|^2 d\x \leq
\left(\frac{C_M}{T-2\theta H}  +
\beta \right) \int_{0}^{T}
 \|\partial_{t}p(\x,t)\|^2_{L^2(\Gamma_m)} dt\\
  + \int_{0}^{T} \| \partial_x p(\x,t)\|^2_{L^2(\Gamma_m)} dt,
\eea
and
\bea
\int_\Omega c^{-2}(y) | f_1(\x)|^2 d\x \leq
\left(\frac{C_M}{T-2\theta H}  +
\beta \right) \int_{0}^{T}
 \|\partial_{t}p_k(\x,t)\|^2_{L^2(\Gamma_m)} dt\\
  + \int_{0}^{T} \| \partial_x p(\x,t)\|^2_{L^2(\Gamma_m)} dt,
\eea
with 
\bea
C_M = H e^{\int_0^H  c^{2}(s)|\partial_y (c^{-2}(s))|ds } 
(c^{-2}(H) +\beta^2).
\eea

\end{theorem}
\proof
The estimates are direct consequences of  Proposition~\ref{pr1} and Proposition~\ref{pr2}. 
The fact that  the Fourier series of $p(\x, t)$  has a finite number of terms justifies the
regularity of the solution $p(\x, t)$, and  allow interchanging the order between the Fourier 
series and the  integral over $(0,T)$. 
\endproof

Using microlocal analysis techniques it is known that the boundary observability   
in a   rectangle  holds  if the set of boundary observation necessarily
contains at  least two adjacent sides \cite{BaLeRa-SIAM92,Burq-BSMF98}. Then, we expect that the 
the Lipschitz stability  estimate in Theorem  \ref{mainacoustic} will deteriorate
when the number of modes $N$ becomes larger. In fact the series on the right
side does  not converge because  $\partial_xp(\x,t)$ does not belong in general
to $L^2(\Gamma_m\times (0,T))$.   We here  provide a h\"older stability estimate 
that corresponds to the boundary observability on only one side of the rectangle. 

\begin{theorem} \label{observability}
Assume that $c(y) \in W^{1,\infty}(0,1)$ with $ 0<c_m\leq c^{-2} (y)$, and
$(f_0, f_1)\in H^2(\Omega)\times H^1(\Omega)$ satisfying 
$\|f_0\|_{H^1}, \|f_1\|_{H^2} \leq \tilde M$. 
Let  $\theta = \sqrt{\|c^{-2}\|_{L^\infty}}$  and $T> 2\theta H$. Then
\bea
\int_\Omega  |\nabla f_0(\x)|^2 d\x \leq
\left(\frac{C_M}{T-2\theta H}  +
\beta \right) \int_{0}^{T}
 \|\partial_{t}p(\x,t)\|^2_{L^2(\Gamma_m)} dt\\
  + C_{\tilde M}\left( \int_{0}^{T} \| p(\x,t)\|^2_{H^{\frac{1}{2}}
  (\Gamma_m)} dt\right)^{\frac{2}{3}},
\eea
and
\bea
\int_\Omega c^{-2}(y) | f_1(\x)|^2 d\x \leq
\left(\frac{C_M}{T-2\theta H}  +
\beta \right) \int_{0}^{T}
 \|\partial_{t}p_k(\x,t)\|^2_{L^2(\Gamma_m)} dt\\
  + C_{\tilde M}\theta^{\frac{2}{3}}
  \left( \int_{0}^{T} \| p(\x,t)\|^2_{H^{\frac{1}{2}}(\Gamma_m)} dt\right)^{\frac{2}{3}},
\eea
with 
\bea
C_M = H e^{\int_0^H  c^{2}(s)|\partial_y (c^{-2}(s))|ds } 
(c^{-2}(H) +\beta^2), \;
C_{\tilde M}= 2\tilde M^{\frac{2}{3}}.
\eea
\end{theorem}
\proof
The proof is again based on the results of Proposition \ref{pr2}. We first
deduce from Proposition \ref{preg} that $\partial_t p(\x,t)\in L^2(\Gamma_m)$.
Now, define 
\bea
 f_j^N(\x)&=& \sum_{|k| \leq N} f_{jk}(y) \varphi_k(x) \qquad 
 \x\in \Omega\qquad j=0, 1,
\eea
with $f_{jk}(y)$ are the Fourier coefficients  of $f_j(\x)$, and 
 $N$ being a large positive integer.\\
 
Consequently
 \bean 
\int_\Omega  |\nabla f_0(\x)|^2 d\x \leq  \int_\Omega  |\nabla f_0^N(\x)|^2 d\x
+\frac{\tilde M^2}{\lambda_{N+1}^2},\\
\int_\Omega  c^{-2}(y) | f_1(\x)|^2 d\x \leq  \int_\Omega c^{-2}(y) | f_1^N(\x)|^2  d\x
+\frac{\theta^2 \tilde M^2}{\lambda_{N+1}^2},
 \eean
 for $N$ large. Applying now  Proposition \ref{pr2}  to $(f_0^N, f_1^N)$, gives
 \bea 
\int_\Omega  |\nabla f_0(\x)|^2 d\x \leq  
\left(\frac{C_M}{T-2\theta H}  +
\beta \right) \int_{0}^{T}
 \|\partial_{t}p(\x,t)\|^2_{L^2(\Gamma_m)} dt\\
  + \lambda_N \int_{0}^{T} \| p(\x,t)\|^2_{H^{\frac{1}{2}}(\Gamma_m)} dt
+\frac{\tilde M^2}{\lambda_{N+1}^2},\\
\int_\Omega  c^{-2}(y) | f_1(\x)|^2 d\x \leq 
\left(\frac{C_M}{T-2\theta H}  +
\beta \right) \int_{0}^{T}
 \|\partial_{t}p(\x,t)\|^2_{L^2(\Gamma_m)} dt\\
  + \lambda_N \int_{0}^{T} \|  p(\x,t)\|^2_{H^{\frac{1}{2}}(\Gamma_m)} dt
+\frac{\theta^2 \tilde M^2}{\lambda_{N+1}^2}.
 \eea
 
 By minimizing the right hand terms  with respect to the value  of $\lambda_N$, we obtain the 
 desired results.

\endproof

\section{The optical inversion} \label{opticalinversion}

Once the initial pressure $f_0(\x)$, generated by the optical wave has been reconstructed, a second step consists of  determining the optical properties in the sample. Although this second step has not been well studied in biomedical literature due its complexity, it is of importance in applications. In fact  the optical parameters  are very sensitive to the tissue condition and their values for healthy and unhealthy tissues are extremely different.
 
The second inversion is to determine the coefficients $(D(y), \mu_a(y))$ from 
the initial  pressures recovered in the first inversion, that is,
${ h}_j(\x) = \mu_a(y) u_j(\x),  \; \x\in \Omega, \; j=1, 2$.\\

For simplicity, we will  consider $g_j(\x) = \varphi_{k_j}(x), \; j=1,2$
with $k_1$ and $k_2$ are two distinct Fourier eigenvalues 
that  are large enough. We specify
how large they should be later in the analysis.

The main result of this section is the following. 
\begin{theorem} \label{mainopticinversion}
Let $(D,\mu_a)$,\,  $(\widetilde D, \tilde \mu_a)$
in $\mathcal O_M$, and $k_i, \, i=1,2$ be two distinct integers.
 Denote $u_{k_i}, \, i=1,2$
 and $\tilde u_{k_i}, \, i=1,2 $  the solutions  to 
 the system~\eqref{equk} for $g_{i} = \varphi_{k_i}, i=1,2$, with coefficients 
 $(D,\mu_a)$ and  $(\widetilde D, \tilde \mu_a)$ respectively.
 Assume that $D(H)= 
 \widetilde D(H)$,  $D^\prime(H)= 
 \widetilde D^\prime(H)$, $\mu_a^\prime(H)= 
 \tilde \mu_a^\prime(H),\,$  $k_1<k_2$, and $k_1$ is large enough. 
  Then, there exists
a constant $C>0$ that only depends on 
$(\mu_0, D_0, k_1, k_2, M, L, H)$,  such that 
 the following stability estimates  hold.
 \bea
\| \underline u_{m}^2(D- \tilde D) \|_{C^0} \leq C
\left(\|h_1- \tilde h_1\|_{C^1}+ 
\|h_2- \tilde h_2\|_{C^1}\right),\\
\| \underline u_{m}^2(\mu_a- \tilde \mu_a) \|_{C^0} \leq C
\left(\|h_1- \tilde h_1\|_{C^1}+ 
\|h_2- \tilde h_2\|_{C^1}\right).
\eea
 \end{theorem}

Classical elliptic operator theory implies the following result
for the direct problem~\cite{McLean-Book00}.
\begin{proposition}
Assume $(D,\mu_a)$  be  in $\mathcal O_M$ and 
$g\in V_{\Gamma_m}$. 
Then, there exists a unique solution $u\in V$ to the 
system~\eqref{maineqO}.
It verifies
 \bea
 \| u\|_{H^1(\Omega)} \leq C_0 \|g\|_{H^{\frac{1}{2}}(\Gamma_m)}^2,
 \eea
 where $C_0= C_0(\mu_0, D_0, M,  L, H)>0$.
\end{proposition}

For  $g(\x) = \varphi_k(x)$, the unique solution $u$
 has the following decomposition
 \bea
u(\x)&=& u_k(y) \varphi_k(x) \qquad 
 \x \in \Omega,
\eea
where  $u_k(y)$ satisfies the following one dimensional elliptic equation
\bean \label{equk}
\left \{ \ba{llllccc} 
-\left(D(y)u^\prime(y)\right)^\prime +(\mu_a(y)+
\lambda_k^2 D(y) )u(y) = 0  &y\in (0,H), \\
u(H)= 1, \;\;
 u(0) =0,&
 \ea
\right.
\eean
Next we will derive
 some useful properties of the solution to the system ~\eqref{equk}.

 \begin{lemma} \label{uinf}Let $u(y)$ be the unique solution to 
 the system ~\eqref{equk}. Then $u(y) \in C^2([0, H])$
 and there exists a constant $b= b(\mu_0, D_0, M,  L, H)>0$ such that
  $\|u\|_{C^2} \leq b$ for all $(D,\mu_a)\in \mathcal O_M$. In addition
 the following inequalities hold for $k$ large enough.  
 \bea 
  \underline u_m(y) \leq  u(y)\leq
 \overline u_M(y),
 \eea
 for $0\leq y\leq H $, where 
 \bea
 \underline u_m(y) = \frac{D^{\frac{1}{2}}(H)}{D^{\frac{1}{2}}(y)}
 \frac{\sinh(\kappa_m^{\frac{1}{2}} y)}
 { \sinh(\kappa_m^{\frac{1}{2}}H)},\quad 
\overline u_M(y) = 
 \frac{D^{\frac{1}{2}}(H)}{D^{\frac{1}{2}}(y)} \frac{\sinh(\kappa_M^{\frac{1}{2}} y)}
 { \sinh(\kappa_M^{\frac{1}{2}} H)},\\ 
 \kappa_m = \min_{0\leq y\leq H}
 \left( \frac{(D^{\frac{1}{2}})^{\prime\prime}}{D^{\frac{1}{2}}}+ \frac{\mu_a}{D}+
\lambda_k^2\right),\quad 
 \kappa_M = \max_{0\leq y\leq H} \left( 
 \frac{(D^{\frac{1}{2}})^{\prime\prime}}{D^{\frac{1}{2}}}+ \frac{\mu_a}{D}+
\lambda_k^2\right).
\eea
 \end{lemma} \proof
We first make 
 the Liouville change of variables and introduce the function 
 \[
 v(y) = \frac{D^{\frac{1}{2}}(y)}{D^{\frac{1}{2}}(H)} 
 u(y).
 \]
 Forward
 calculations show that $v(y)$ is the unique solution to the following system.
 \bean \label{eqvk}
\left \{ \ba{llllccc} 
-v^{\prime\prime}(y)+ \kappa(y)v(y) = 0  &y\in (0,H), \\
v(H)= 1, \;\;
v(0) =0,&
 \ea
\right.
\eean
where $$\kappa(y) = \frac{(\sqrt D)^{\prime\prime}}{\sqrt D}+ \frac{\mu_a}{D}+
\lambda_k^2.  $$
Assume now that $k$ is large enough such that
 $\kappa_m>0$, and
let $\underline v_m(y)$ and $\overline v_M(y)$ be the solutions to the system~\eqref{eqvk}
when we replace  $\kappa(y)$ by respectively the constants $\kappa_m$
 and $\kappa_M$.
They are {explicitly} given by 
\bea
\underline v_m(y) &=&\frac{\sinh(\sqrt{\kappa_m} y)}
 { \sinh(\sqrt{\kappa_m} H)},\\
 \overline v_M(y) &=&\frac{\sinh(\sqrt{\kappa_M} y)}
 { \sinh(\sqrt{\kappa_M} H)}.
\eea
The maximum principle~\cite{McLean-Book00} implies that $0 <v(y), \underline
v_m(y), 
\overline v_M(y)<1$
for  $0< y< H $.\\

By applying again the maximum principle on the differences 
$v-\underline v_m$ and $v-\overline v_M$ we deduce that $ \underline 
v_m(y)<v(y)<\overline v_M(y)$
for  $0\leq  y \leq H $, which
leads to the desired lower and upper bounds.\\  

We deduce from the regularity of the coefficients $D$ and $\mu_a$  and
the {classical} elliptic regularity \cite{McLean-Book00} that $u\in H^{3}(0,
H)$. Moreover there exist a constant $b>$ that only depends on 
$(\mu_0, D_0, M,  L, H)$ 
such that
\bean \label{boundsonu}
 \|u\|_{H^{3}} \leq b.
\eean
Consequently the uniform  $C^2$ bound of $u$ can be obtained 
using the  continuous Sobolev embedding of   $H^{3}(0,
H)$  into $C^2([0,H])$ \cite{AdFo-Book03}. 

\endproof
\begin{lemma}\label{uprimeinf} Let $(D,\mu_a)\in \mathcal O_M$,
and $u(y)$ be the unique solution to 
 the system~\eqref{equk}. Then,
 for $k$ large enough there exists a constant
 $\varrho = \varrho(D_0, \mu_0, M, k) >0$ such that 
 \bea
 u^\prime(y) \geq \varrho,
 \eea
 for  $0\leq y \leq H$.
\end{lemma}
\proof
Since $0$ is the global minimum of
$u(y)$, we have $u^\prime(0)>0$.  Moreover for  $k$
large enough, Lemma~\ref{uinf} implies
that 
 \bea 
u(y) \geq \frac{\sqrt{D(H)}}{\sqrt{D(y)}} \frac{\sinh(\sqrt{\kappa_m} y)}
 { \sinh(\sqrt{\kappa_m} H)},
  \eea
for all $y \in [0, H]$. Therefore 
\bea
u^\prime(0) \geq \frac{\sqrt{D(H)}}{\sqrt{\|D\|_{L^\infty }}} \frac{\sqrt{\kappa_m}}
 { \sinh(\sqrt{\kappa_m} H)}
\eea

Now integrating equation~\eqref{equk} over $(0, y)$
we obtain 
\bea
D(y)u^\prime(y) =  D(0)u^\prime(0) +\int_0^y
(\mu_a(s)+ \lambda_k^2 D(s))u(s)ds  
\eea
\bea
D(y)u^\prime(y) \geq   D(0)u^\prime(0) +\int_0^y
(\mu_a(s)+ \lambda_k^2 D(s)) \frac{\sqrt{D(H)}}{\sqrt{D(s)}} 
\frac{\sinh(\sqrt{\kappa_m} s)}
 { \sinh(\sqrt{\kappa_m} H)}ds  \\
\geq   \frac{\sqrt{D(H)}}{\sqrt{\|D\|_{L^\infty }}} \frac{\sqrt{\kappa_m} D_0}
 { \sinh(\sqrt{\kappa_m} H)}
+
(\mu_0+ \lambda_k^2 D_0) \frac{\sqrt{D(H)}}{\sqrt{ \|D\|_{L^\infty }} }
\frac{\cosh(\sqrt{\kappa_m}y) -1}
 { \sinh(\sqrt{\kappa_m} H)}.
\eea
Taking into account the  explicit expression of $\kappa_m$ finishes the proof.
\endproof

 Since the illumination are chosen to coincide 
 with the Fourier basis functions $\varphi_{k_j}, \, j=1,2$, 
 the data
 ${\bf h}_j(\x), \, j=1,2$,  can be rewritten as $ {\bf h}_j(\x)
  = h_j(y)\varphi_{k_j}(x), \, j=1,2$, where $h_j(y) =  \mu_a(y)
  u_{k_j}(y)$.  \\
  
 Therefore the optical inversion is reduced to
 the problem of identifying  the optical pair    
 $(D,\mu_a)$ from the knowledge of the pair
 $(h_1(y), h_2(y))$ over $(0,H)$. \\

  Let $(D,\mu_a)$,\,  $(\widetilde D, \tilde \mu_a)$ be two
 different pairs in $\mathcal O_M$, and denote $u_{k}$
 and $\tilde u_{k}$  the solutions  to 
 the system~\eqref{equk}, with coefficients 
 $(D,\mu_a)$ and  $(\widetilde D, \tilde \mu_a)$ respectively. \\
 
 We deduce from Lemma~\ref{uinf} that $\frac{1}{u_k}$ 
and $\frac{1}{\tilde u_k}$ lie in $L^p(0, H)$ for $0<p<1$. 
Unfortunately
for or $0<p<1$ the usual $\|\cdot\|_{L^p}$ is not 
anymore a norm  on the  vector  space  
$L^p(0, H)$ because it does not satisfy the triangle inequality
  (see for instance \cite{AdFo-Book03}). In contrast with triangle inequality
  H\"older inequality holds for $0<p<1$, and we have 
\bean \label{holderin}
\| \frac{v}{u_k} \|_{L^r} \leq \| \frac{1}{u_k} \|_{L^p} 
\|v \|_{L^q}, 
\eean
for all $v \in L^q(0, H)$ with $\frac{1}{r} = \frac{1}{p}+\frac{1}{q}$.\\

Consequently $h =  \frac{h_2}{h_1}= 
\frac{u_{k_2}}{u_{k_1}}$ 
 can be considered as a distribution that coincides with
 a $C^2$ function over $(0,H)$. 
A forward calculation shows that  $h$ satisfies the equation 
\bean \label{heq}
-\left(Du_{k_1}^2 h^\prime \right)^\prime + Du_{k_1}^2 h
(\lambda_2^2-\lambda_1^2)= 0,
\eean
over $(0, H)$. \\

Since $u_{k_j}, \, j=1,2$ are in $C^2([0, H])$, an asymptotic
analysis of $ D(y)u_{k_1}^2(y) h^\prime(y)$ at $0$ and 
the results of
Lemma \ref{uinf}, gives
\bea
\lim_{y \to 0}Du_{k_1}^2 h^\prime=0. 
\eea
Similarly, we have
\bea
\lim_{y \to 1}h = 1. 
\eea

Integrating the equation \eqref{heq} over $(0, y)$, we get
\bea
D(y)u_{k_1}^2(y) h^\prime(y) = (\lambda_2^2-\lambda_1^2)\int_0^y D(s)
 u_{k_1}(s)u_{k_2}(s) ds.
\eea
Dividing both sides by $D(y)u_{k_1}^2(y) $, and using again Lemma \ref{uinf}, imply
\bea
 h^\prime(y) \geq  (\lambda_2^2-\lambda_1^2) D^{-1}(y)\overline{u}_{M}^{-2}(y)\int_0^y D(s)
 \underline{u}_{m}^{2}(s) ds,
\eea
which leads to 
\bean \label{hinf}
h^\prime(y) \geq  (\lambda_2^2-\lambda_1^2)
\min_{y\in (0, H)} D^{-1}(y)\overline{u}_{M}^{-2}(y)\int_0^y D(s)
 \underline{u}_{m}^{2}(s) ds>0. 
\eean

The right hand constant is strictly positive and only depends on $D_0, \mu_0, M, H, L$ 
and $k$.\\

Now back to the optical inversion. The equation \eqref{heq}
can be written as
\bea
-(Du_{k_1}^2)^\prime h^\prime  + \left( h
(\lambda_2^2-\lambda_1^2) - h^{\prime\prime}\right)
Du_{k_1}^2= 0,
\eea
over $(0, H)$.  Dividing both sides by $Du_{k_1}^2 h^\prime$,
and integrating over $(0, y)$,
we obtain
\bean \label{feq}
D(y)u_{k_1}^2(y)= h(0) -h(y)+e^{(\lambda_2^2-\lambda_1^2)
\int_0^y \frac{h}{h^\prime} ds}.
\eean
{This allows us to show the following result.}
\begin{lemma} \label{stabinter} Under the assumptions of Theorem 
\ref{mainopticinversion}, 
there exists
 a constant $C= C(\mu_0, D_0, k_1, k_2, M, L, H)>0$ such that 
 the following inequality holds.
 \bea
  \| D  u_{k_1}^2-\widetilde D
 \tilde u_{k_1}^2 \|_{C^0} \leq  C \left(\|h_1-\tilde h_1 \|_{C^1} +
  \|h_2-\tilde h_2 \|_{C^1}\right).
 \eea

 \end{lemma}
\proof 
Recall that the relation \eqref{feq} is also
valid for the pair  $(\widetilde D, \tilde \mu_a)$, that is
\bea 
\widetilde D(y)\tilde u_{k_1}^2(y)= \tilde h(0) -\tilde h(y)+
e^{(\lambda_2^2-\lambda_1^2)
\int_0^y \frac{\tilde h}{\tilde h^\prime} ds},
\eea
where $\tilde h =  \frac{\tilde u_{k_2}}{\tilde u_{k_1}}$.
Taking the difference between the last equation and 
the equation \eqref{feq} we find 
\bea 
\| D  u_{k_1}^2-\widetilde D
 \tilde u_{k_1}^2 \|_{C^0} \leq \| h-\tilde h\|_{C^0} 
+ (\lambda_2^2-\lambda_1^2)H \left\| \frac{h}{h^\prime}
- \frac{\tilde h}{\tilde h^\prime} \right\|_{C^0}
e^{(\lambda_2^2-\lambda_1^2)H \left(\left\| \frac{h}{h^\prime}
\right\|_{C^0}+\left\|\frac{\tilde h}{\tilde h^\prime} \right\|_{C^0}\right)}.
\eea
We then deduce the result from Lemma \ref{uinf}
and inequality \eqref{hinf}.

\endproof

Now, we are ready to prove the main stability result of this
section.  We remark as in \cite{BaRe-IP11}, that $ 
\frac{1}{u_{k_1}} $ is a solution to the
following equation.
\bea
-\left(Du_{k_1}^2  \frac{1}{u_{k_1}} ^\prime\right)^\prime +
\lambda_{k_1}^2 Du_{k_1}^2 \frac{1}{u_{k_1}}  = h_1,  &y\in (0,H).
\eea
Since $ \frac{1}{\tilde u_{k_1}} $ solves the same type of 
equation, we obtain that  $w = 
\frac{1}{u_{k_1}} - \frac{1}{\tilde u_{k_1}} $, is the
solution to the following system
\bean
\left \{
\ba{llcc}
-\left(Du_{k_1}^2 w ^\prime\right)^\prime +
\lambda_{k_1}^2 Du_{k_1}^2 w  =  e,  &y\in (0,H),\\
w(H) = 0,\; w^\prime(H) = \frac{1}{\mu_a(H)}(\tilde h_1^\prime(H) -
h_1^\prime(H)) ,  &yw(y) \in L^2(0, H),
\ea
\right.
\eean
where 
\bea
e= -\left((Du_{k_1}^2- \widetilde D\tilde u_{k_1}^2 )
\frac{1}{\tilde u_{k_1}} ^\prime\right)^\prime
+\lambda_{k_1}^2(Du_{k_1}^2- \widetilde D\tilde u_{k_1}^2)
 \frac{1}{\tilde u_{k_1}} +h_1- \tilde h_1.
\eea
We  remark that to solve this system  
we have to deal with two main difficulties, the  first is that the
operator is elliptic degenerate,  and the second
is that the solution $w(y)$ may be unbounded at $y=0$. \\

Multiplying by $\mathrm{sign}(w)$, and integrating over $(s, H)$ 
the first equation of the
system leads to
\bea
D(s)u_{k_1}^2(s) |w|^\prime(s)= \\
\mathrm{sign}(w)\left(\frac{D(H)}{\mu_a(H)}
(\tilde h_1^\prime(H)-h_1^\prime(H))
+ \int_s^H e(y) dy -
\lambda_{k_1}^2 \int_s^H  D(y)u_{k_1}^2(y) w(y) dy\right).
\eea
Integrating again over $(t, H)$ gives
\bea
\int_t^H D(s)u_{k_1}^2(s) |w|^\prime(s) ds \leq \\
 \frac{MH}{\mu_0} \|h_1- \tilde h_1\|_{C^1}
+  \int_0^H \left|\int_s^H e(y) dy\right| ds+
\lambda_{k_1}^2 H \int_0^H  D(y)u_{k_1}^2(y) |w|(y) dy.
\eea
Since $u_{k_1}^\prime>0$ over $(0, H)$ (Lemma \ref{uprimeinf}), $u_{k_1}$
is increasing, and we have
\bea
 D(t)u_{k_1}^2(t) |w|(t) \leq \\ 
 \frac{MH}{\mu_0}  \|h_1- \tilde h_1\|_{C^1}
+  \int_0^H\left| \int_s^H e(y) dy\right| ds  +
\lambda_{k_1}^2 H \int_0^H  D(y)u_{k_1}^2(y) |w|(y) dy.
\eea

Now, we focus on the second term on the right hand side. 
\bean\label{o1}
\int_0^H \left|\int_s^H e(y) dy\right| ds  \leq \\
\int_0^H |Du_{k_1}^2- \widetilde D\tilde u_{k_1}^2 |
  \frac{|\tilde u_{k_1}^\prime| }{|\tilde u_{k_1}^2|} ds + 
\lambda_{k_1}^2 H\int_0^H |Du_{k_1}^2- 
\widetilde D\tilde u_{k_1}^2|
 \frac{1}{|\tilde u_{k_1}|} dy+  H \|h_1- \tilde h_1\|_{C^0}.\nonumber
 \eean
Using the estimates in  Lemma \ref{stabinter},  we find

\bean\label{o2}
\int_0^H \left|\int_s^H e(y) dy\right| ds  \leq 
C \left(\|h_1-\tilde h_1 \|_{C^1} +  \|h_2-\tilde h_2 \|_{C^1}\right).
\eean

Combining inequalities \eqref{o1} and \eqref{o2},  leads to
\bean \label{wf}
D(t)u_{k_1}^2(t)|w|(t)\leq\\ C_1\left(\|h_1- \tilde h_1\|_{C^1}+ \|h_2- 
\tilde h_2\|_{C^1}\right) +   C_2\int_t^H 
D(y)u_{k_1}^2(y)|w|(y) dy,\nonumber
\eean
for $0\leq t \leq H$.\\

Using Gronwall's inequality we get 
\bea
D(t)u_{k_1}^2(t)|w|(t)\leq C_1 e^{C_2 \int_0^H 
D(y)u_{k_1}^2(y)  dy}
 \left(\|h_1- \tilde h_1\|_{C^1}+ \|h_2- 
\tilde h_2\|_{C^1}\right),
\eea
for $0\leq t \leq H$.\\
Finally, we obtain the following estimate. 
\bean \label{tweq}
{u}_{k_1}^2(t)|w|(t)\leq D_0^{-1}C_1 e^{C_2 \int_0^H 
D(y)u_{k_1}^2(y)  dy}
 \left(\|h_1- \tilde h_1\|_{C^1}+ \|h_2- 
\tilde h_2\|_{C^1}\right),
\eean

The following Lemma is a direct consequence
of the  previous inequality. 

\begin{lemma} \label{stabinter2}
Under the assumptions of Theorem \ref{mainopticinversion}, there exists
 a constant $C = C(\theta, \mu_0, D_0, k_1, k_2, M, L, H)>0$
 such that  the following inequality holds.
 \bea
  \| {u}_{k_1}(u_{k_1}-
 \tilde u_{k_1})\|_{C^0} \leq 
  C \left(\|h_1- \tilde h_1\|_{C^1}+ \|h_2- 
\tilde h_2\|_{C^1}\right).
 \eea
  \end{lemma}

\proof(Theorem \ref{mainopticinversion})
Recall that $h_1= \mu_a u_{k_1}$ and
$\tilde h_1= \tilde \mu_a \tilde u_{k_1}$ over $(0, H)$. \\

Therefore 
\bea
u_{k_1}^2|\mu_a - \tilde \mu_a| \leq  u_{k_1}| h_1-\tilde h_1|+ 
\tilde \mu_au_{k_1} |u_{k_1}-\tilde u_{k_1}|.
\eea
Lemma \ref{stabinter2} implies 

\bean
  \| u_{k_1}^2(\mu_a - \tilde \mu_a) \|_{C^0} \leq 
  C \left(\|h_1- \tilde h_1\|_{C^1}+ \|h_2- 
\tilde h_2\|_{C^1}\right).
 \eean

A simple calculation yields 

\bea
u_{k_1}^2|D-\widetilde D| \leq \widetilde D
|u_{k_1}^2-\tilde u_{k_1}^2| + |Du_{k_1}^2-\widetilde D\tilde u_{k_1}^2|, 
\eea
over $(0, H)$.

 Lemma \ref{stabinter} 
and \ref{stabinter2} leads to 
\bean \label{ineqDD}
\| u_{k_1}^2(D- \tilde D) \|_{C^0} \leq C
\left(\|h_1- \tilde h_1\|_{C^1}+ 
\|h_2- \tilde h_2\|_{C^1}\right).
\eean

Applying the bounds in Lemma \ref{uinf}, we obtain the wanted results.
\endproof

\section{Proof of Theorem~\ref{mainglobalinversion}}

The main idea here is to {combine} the stability results of the acoustic
and optic inversions in a result that shows how the reconstruction of
the optical coefficients is sensitive to the noise in the measurements
of the acoustic waves.\\

The {principal} difficulty is that the vector spaces used in both stability estimates are not 
the same due to the difference in the techniques used to derive them. We will  {use} interpolation
 inequality between Sobolev spaces to overcome this difficulty.

We  deduce from  the uniform bound on the solutions $u_i, \, i=1,2$ (see for instance \eqref{boundsonu} in the proof of Lemma~\ref{uinf}) that 
\bean \label{hboundsb}
\|h_i\|_{H^{3}}, \|\tilde h_i\|_{H^{3}} \leq Mb, \quad i=1,2,
\eean
for all pairs $(D,\mu_a)$and  $(\widetilde D,\mu_a)$ in $\mathcal
O_M$.

The Sobolev interpolation inequalities and  embedding theorems~\cite{AdFo-Book03} imply 
\bea
\|h_i-\tilde h_i\|_{C^1} \leq C \|h_i-\tilde h_i\|_{H^{2}}
\leq \widetilde C   \|h_i-\tilde h_i\|_{H^{1}}^{\frac{1}{2}}
\|h_i-\tilde h_i\|_{H^{3}}^{\frac{1}{2}},\quad i=1,2,
\eea
which combined with \eqref{hboundsb} gives
\bean \label{link}
\|h_i-\tilde h_i\|_{C^1} \leq \tilde{\tilde{C}}   
\|h_i-\tilde h_i\|_{H^{1}}^{\frac{1}{2}}, \quad i=1,2. 
\eean

Since the acoustic inversion is linear we obtain from 
Theorem~\ref{mainacoustic} (or Proposition~\ref{pr2}) that,
\bea
\lambda_{k_i}^2\int_0^H  | h_i-\tilde h_i|^2 dy \leq
\left(\frac{C_M}{T-2\theta H}  +
\beta \right) \int_{0}^{T}
 |\partial_{t}p_{i}(H,t)- \partial_{t}\tilde p_{i}(H,t)|^2 dt\\  +
 \lambda_{k_i}^2\int_{0}^{T} |p_i(H,t) - \tilde p_i(H,t) |^2 dt,
\eea
for $i=1,2$, and 
\bea
\int_0^H  c^{-2}(y) |h_i^\prime- \tilde h_i^\prime|^2 dy   \leq
 \left(\frac{C_M}{T-2\theta H}  +\beta \right) \int_{0}^{T}
 |\partial_{t}p_{i}(H,t)- \partial_{t}\tilde p_{i}(H,t)|^2 dt \\ +
 \lambda_{k_i}^2\int_{0}^{T} |p_i(H,t) - \tilde p_i(H,t) |^2 dt,
\eea
for $i=1,2$.\\

Consequently,
\bea
\|h_i-\tilde h_i\|_{C^1} \leq \\ \tilde{C} \left( \int_{0}^{T}
 \left(\frac{C_M}{T-2\theta H}  +\beta \right)|\partial_{t}p_{i}(H,t)- 
\partial_{t}\tilde p_{i}(H,t)|^2 +
 \lambda_{k_i}^2 |p_i(H,t) - \tilde p_i(H,t) |^2
 dt\right)^{\frac{1}{4}},
\eea
for $i=1,2$.

Using the optical stability estimates in Theorem~\ref{mainopticinversion},  we obtain
\bea
\| \underline u_{m}^2(\mu_a- \tilde \mu_a) \|_{C^0}  \leq \\ 
  \tilde{\tilde{C}}
\left( \sum_{i=1}^2 \int_{0}^{T}
 \left(
\frac{C_M}{T-2\theta H}  +\beta \right)|\partial_{t}p_{i}(H,t)- 
\partial_{t}\tilde p_{i}(H,t)|^2 +
 \lambda_{k_i}^2 |p_i(H,t) - \tilde p_i(H,t) |^2
 dt\right)^{\frac{1}{4}},
\eea
 and 
\bea
\| \underline u_{m}^2(D- \tilde D) \|_{C^0} \leq \\ 
\tilde{\tilde{C}}
\left( \sum_{i=1}^2\int_{0}^{T}
 \left(\frac{C_M}{T-2\theta H}  +\beta \right)|\partial_{t}p_{i}(H,t)- 
\partial_{t}\tilde p_{i}(H,t)|^2 +
 \lambda_{k_i}^2 |p_i(H,t) - \tilde p_i(H,t) |^2
 dt\right)^{\frac{1}{4}},
\eea
which ends the proof.
\endproof

\section{Proof of Proposition~\ref{pr1}}

Multiplying the first equation of the system~\eqref{mainpy} 
by $y\partial_{y}p(y,t)$ 
and integrating by part one time over $(0,T)$, we obtain 
\bea
\int_0^T  |\partial_{y}p(H,t)|^2 dt
 = \int_0^T \int_0^H |\partial_{y}p(y,t)|^2 dy dt-
 2 \int_0^T \int_0^H c^{-2} \partial_{tt}p(y,t) y \partial_{y}p(y,t) dy dt\\
- 2\lambda_k^2 \int_0^T \int_0^H p(y,t) y \partial_{y}p(y,t) dy dt= A_1+A_2+A_3.
\eea
In the rest of the proof we shall derive bounds
of  each of the constants $A_i, i=1,2,3, $ in terms of the energy $E(0)$.
Due to the energy decay~\eqref{energderivative}, we have
\bea
|A_1|\leq TE(0).
\eea
Integrating by part again over $(0,T)$ in the integral $A_2$, we get
\bea 
A_2= -\int_0^T \int_0^H yc^{-2} \partial_{y}|\partial_{t}p(y,t)|^2  dy dt\\
+2\int_0^H yc^{-2} \partial_{t}p(y,T) \partial_{y}p(y,T) dy
-2\int_0^H yc^{-2} \partial_{t}p(y,0) \partial_{y}p(y,0) dy.
\eea
Integrating by part now over  $(0,H)$, we find
\bea 
A_2 + H c^{-2}(H) \int_0^T |\partial_{t}p(H,t)|^2 dt =
  \int_0^T \int_0^H \partial_{y}(yc^{-2}) |\partial_{t}p(y,t)|^2  dydt\\
+2\int_0^H yc^{-2} \partial_{t}p(y,T) \partial_{y}p(y,T) dy
-2\int_0^H yc^{-2} \partial_{t}p(y,0)  \partial_{y}p(y,0) dy,
\eea
which leads to the following inequality
\bea 
\left |A_2 +H c^{-2}(H) \int_0^T |\partial_{t}p(H,t)|^2 dt \right |
 \leq\\ \left\|c^2\partial_{y}(yc^{-2}(y)) \right \|_{L^\infty} TE(0) 
 +H \|c^{-1}(y)\| (E(T) + E(0)).
\eea
Using the energy decay~\eqref{energderivative}, we finally obtain
\bea
\left |A_2 +H c^{-2}(H) \int_0^T |\partial_{t}p(H,t)|^2 dt \right |
\leq\\ \left( (1+(1+\frac{H}{c_m})\| c^{-2}\|_{ W^{1,\infty}})T+ 2H
\| c^{-2}\|_{ L^{\infty}}^{1/2} \right) E(0).
\eea
Similar arguments  dealing with the integral $A_3$ 
show that
\bea
|A_3|\leq H\lambda_kT E(0). 
\eea
Combining all the previous estimates on the constants 
$A_i, i=1,2,3, $ achieve the proof.

\section{Proof of Proposition~\ref{pr2}}

Let  $\theta = \sqrt{\|c^{-2}\|_{L^\infty}}$  and $T> 2\theta H, $ and 
introduce the following function
\bea
\Phi(y) \hspace{-3mm}&=& \hspace{-3mm} 
\int_{\theta y}^{T-\theta y}\left(c^{-2}(H-y)  |\partial_{t}p(H-y,t)|^2  +
 |\partial_{y}p(H-y,t)|^2 +\lambda_k^2 |p(H-y,t)|^2 \right) dt,\\
 &=&\int_{\theta y}^{T-\theta y} \varphi(y,t) dt,
\eea
for $0\leq y \leq H$.
We remark that 
\bean \label{phi0}
\Phi(0) &= & (c^{-2}(H) +\beta^2)\int_{0}^{T}
 |\partial_{t}p(H,t)|^2 dt + \lambda_k^2\int_{0}^{T} |p(H,t)|^2 dt.
\eean
On the other hand a forward calculation  of the derivative
of $\Phi(y)$
gives
\bea
\Phi^\prime(y) = \int_{\theta y}^{T-\theta y} \partial_y \varphi(y,t) dt
- \theta \varphi(y, T-\theta y) - \theta \varphi(y, \theta y).
\eea
Integrating by parts in the integral  we deduce that
\bea
\Phi^\prime(y) =   B_\theta(y) 
+ \partial_y (c^{-2}(H-y)) \int_{\theta y}^{T-\theta y} |\partial_{t}p(H-y,t)|^2 dt,
 \eea
 where 
 \bea
  B_\theta(y)  = \left(-2c^{-2}(H-y)\partial_{t}p(H-y,t) 
  \partial_{y}p(H-y,t)\right)\Big|_{t= \theta y}^{t=T-\theta y} \\
 - \theta \left(c^{-2}(H-y)  |\partial_{t}p(H-y,t)|^2  +
 |\partial_{y}p(H-y,t)|^2 +\lambda_k^2 |p(H-y,t)|^2 \right)
 \Big|_{t= \theta y}^{t=T-\theta y}.
 \eea
The choice of $\theta$ implies $B_\theta(y)<0$ for $0\leq y \leq H$. 
Hence, we obtain 
\bea
\Phi^\prime(y) \leq c^{2}(H-y)|\partial_y (c^{-2}(H-y))| 
\int_{\theta y}^{T-\theta y} c^{-2}(H-y) |\partial_{t}p(H-y,t)|^2 dt\\
\leq  c^{2}(H-y)|\partial_y (c^{-2}(H-y))| \Phi(y).
 \eea
 Using Gronwall's inequality we get
\bean \label{ineqphi}
 \Phi(y)\leq e^{\int_0^H  c^{2}(s)|\partial_y (c^{-2}(s))|ds } \Phi(0),
\eean
for $0\leq y \leq H$.\\

We deduce from the energy decay~\eqref{energderivative} that 
\bean  \label{ineqphi1}
(T-2\theta H) E(T) \leq (T-2\theta H) E(T-\theta H)
\leq \int_{\theta H}^{T-\theta H} E(t) dt.
\eean
Rewriting now the right hand side in terms of the
function $\varphi$ we found
\bea
\int_{\theta H}^{T-\theta H} E(t) dt = \int_0^H
\int_{\theta H}^{T-\theta H}  \varphi(y,t) dt dy.
\eea
Since $(\theta H, T- \theta H) \subset (\theta y, T-\theta y)$
for all $0\leq y \leq H$, we have 
\bean  \label{ineqphi2}
\int_{\theta H}^{T-\theta H} E(t) dt \leq 
\int_0^H  \Phi(y) dy.
\eean
Combining inequalities~\eqref{ineqphi}-\eqref{ineqphi1}-\eqref{ineqphi2}, we
 find
\bean \label{fin}
(T-2\theta H) E(T) \leq H e^{\int_0^H  c^{2}(s)|\partial_y (c^{-2}(s))|ds } \Phi(0).
\eean
Back again to the energy derivative~\eqref{energderivative}, and
integrating the equality over $(0,T)$ we obtain
\bea
E(0) = E(T) + \beta \int_0^T  |\partial_{t}p(H,t)|^2 dt. 
\eea
The last equality and energy estimate~\eqref{fin} give
\bea
E(0) \leq (T-2\theta H)^{-1} H e^{\int_0^H  c^{2}(s)|\partial_y (c^{-2}(s))|ds } \Phi(0) +
\beta \int_0^T  |\partial_{t}p(H,t)|^2 dt. 
\eea
Substituting $\Phi(0)$ by its expression in~\eqref{phi0} we finally
find
\bea
E(0) \leq\\ \left((T-2\theta H)^{-1} H e^{\int_0^H  c^{2}(s)|\partial_y (c^{-2}(s))|ds } 
(c^{-2}(H) +\beta^2) +\beta \right) \int_{0}^{T}
 |\partial_{t}p(H,t)|^2 dt\\  + \lambda_k^2\int_{0}^{T} |p(H,t)|^2 dt,
\eea
which combined with the fact that 
\bea
E(0) &=& \int_0^H  \left(c^{-2}(y)  |f_1(y)|^2  +
 |f_0^\prime(y)|^2 +\lambda_k^2 |f_0(y)|^2 \right) dy,
 \eea
finishes the proof. 
\endproof

\section*{Acknowledgments}

The work of KR is partially supported by the US National Science Founda-
tion through grant DMS-1620473. The research of FT was supported in part by
the LabEx PERSYVAL-Lab (ANR-11-LABX- 0025-01). FT would like to thank
the Institute of Computational Engineering and Sciences (ICES) for the provided
support during his visit.

{\small

\begin{thebibliography}{10}

\bibitem{AdFo-Book03}
{\sc R.~A. Adams and J.~F. Fournier}, {\em {Sobolev} {Spaces}}, Academic Press,
  2nd~ed., 2003.

\bibitem{AgKuKu-PIS09}
{\sc M.~Agranovsky, P.~Kuchment, and L.~Kunyansky}, {\em On reconstruction
  formulas and algorithms for the thermoacoustic tomography}, in Photoacoustic
  Imaging and Spectroscopy, L.~V. Wang, ed., CRC Press, 2009, pp.~89--101.

\bibitem{AgQu-JFA96}
{\sc M.~Agranovsky and E.~T. Quinto}, {\em Injectivity sets for the {Radon}
  transform over circles and complete systems of radial functions}, J. Funct.
  Anal., 139 (1996), pp.~383--414.

\bibitem{AmBoCaTaFi-SIAM08}
{\sc H.~Ammari, E.~Bonnetier, Y.~Capdeboscq, M.~Tanter, and M.~Fink}, {\em
  Electrical impedance tomography by elastic deformation}, SIAM J. Appl. Math.,
  68 (2008), pp.~1557--1573.

\bibitem{AmBoJuKa-SIAM10}
{\sc H.~Ammari, E.~Bossy, V.~Jugnon, and H.~Kang}, {\em Mathematical modelling
  in photo-acoustic imaging of small absorbers}, SIAM Rev., 52 (2010),
  pp.~677--695.

\bibitem{AmBrGaJu-SIAM12}
{\sc H.~Ammari, E.~Bretin, J.~Garnier, and V.~Jugnon}, {\em Coherent
  interferometry algorithms for photoacoustic imaging}, SIAM J. Numer. Anal.,
  (2012).

\bibitem{AmBrJuWa-LNM12}
{\sc H.~Ammari, E.~Bretin, V.~Jugnon, and A.~Wahab}, {\em Photo-acoustic
  imaging for attenuating acoustic media}, in Mathematical Modeling in
  Biomedical Imaging II, H.~Ammari, ed., vol.~2035 of Lecture Notes in
  Mathematics, Springer-Verlag, 2012, pp.~53--80.

\bibitem{AmKaKi-JDE12}
{\sc H.~Ammari, H.~Kang, and S.~Kim}, {\em Sharp estimates for {Neumann}
  functions and applications to quantitative photo-acoustic imaging in
  inhomogeneous media}, J. Diff. Eqn., 253 (2012), pp.~41--72.

\bibitem{AmCh-DPDE17}
{\sc K.~Ammari and M.~Choulli}, {\em Logarithmic stability in determining a
  boundary coefficient in an {IBVP} for the wave equation}, Dynamics of PDE, 14
  (2017), pp.~33--45.

\bibitem{AmChTr-PAMS16}
{\sc K.~Ammari, M.~Choulli, and F.~Triki}, {\em Determining the potential in a
  wave equation without a geometric condition. extension to the heat equation},
  Proc. Amer. Math. Soc., 144 (2016), pp.~4381--4392.

\bibitem{AmChTr-arXiv16}
\leavevmode\vrule height 2pt depth -1.6pt width 23pt, {\em H\"older stability
  in determining the potential and the damping coefficient in a wave equation},
  arXiv:1609.06102,  (2016).

\bibitem{Arridge-IP99}
{\sc S.~R. Arridge}, {\em Optical tomography in medical imaging}, Inverse
  Probl., 15 (1999), pp.~R41--R93.

\bibitem{Bal-IO12}
{\sc G.~Bal}, {\em Hybrid inverse problems and internal functionals}, in Inside
  Out: Inverse Problems and Applications, G.~Uhlmann, ed., vol.~60 of
  Mathematical Sciences Research Institute Publications, Cambridge University
  Press, 2012, pp.~325--368.

\bibitem{BaRe-IP11}
{\sc G.~Bal and K.~Ren}, {\em Multi-source quantitative {PAT} in diffusive
  regime}, Inverse Problems, 27 (2011).
\newblock 075003.

\bibitem{BaRe-CM11}
\leavevmode\vrule height 2pt depth -1.6pt width 23pt, {\em Non-uniqueness
  result for a hybrid inverse problem}, in Tomography and Inverse Transport
  Theory, G.~Bal, D.~Finch, P.~Kuchment, J.~Schotland, P.~Stefanov, and
  G.~Uhlmann, eds., vol.~559 of Contemporary Mathematics, Amer. Math. Soc.,
  Providence, RI, 2011, pp.~29--38.

\bibitem{BaUh-IP10}
{\sc G.~Bal and G.~Uhlmann}, {\em Inverse diffusion theory of photoacoustics},
  Inverse Problems, 26 (2010).
\newblock 085010.

\bibitem{BaUh-CPAM13}
\leavevmode\vrule height 2pt depth -1.6pt width 23pt, {\em Reconstructions of
  coefficients in scalar second-order elliptic equations from knowledge of
  their solutions}, Comm. Pure Appl. Math., 66 (2013), pp.~1629--1652.

\bibitem{BaLeRa-SIAM92}
{\sc C.~Bardos, G.~Lebeau, and J.~Rauch}, {\em Sharp sufficient conditions for
  the observation, control and stabilization of waves from the boundary}, SIAM
  J. Cont. Optim., 30 (1992), pp.~1024--1065.

\bibitem{BuMaHaPa-PRE07}
{\sc P.~Burgholzer, G.~J. Matt, M.~Haltmeier, and G.~Paltauf}, {\em Exact and
  approximative imaging methods for photoacoustic tomography using an arbitrary
  detection surface}, Phys. Rev. E, 75 (2007).
\newblock 046706.

\bibitem{Burq-BSMF98}
{\sc N.~Burq}, {\em Contr\^ole de l'\'equation des ondes dans des ouverts
  comportant des coins}, Bulletin de la S.M.F., 126 (1998), pp.~601--637.

\bibitem{CoArBe-IP07}
{\sc B.~T. Cox, S.~R. Arridge, and P.~C. Beard}, {\em Photoacoustic tomography
  with a limited-aperture planar sensor and a reverberant cavity}, Inverse
  Problems, 23 (2007), pp.~S95--S112.

\bibitem{FiHaRa-SIAM07}
{\sc D.~Finch, M.~Haltmeier, and Rakesh}, {\em Inversion of spherical means and
  the wave equation in even dimensions}, SIAM J. Appl. Math., 68 (2007),
  pp.~392--412.

\bibitem{Grisvard-Book85}
{\sc P.~Grisvard}, {\em Elliptic Problems in Nonsmooth Domains}, Pitman
  Publishing Inc., 1985.

\bibitem{Haltmeier-SIAM11}
{\sc M.~Haltmeier}, {\em Inversion formulas for a cylindrical {Radon}
  transform}, SIAM J. Imag. Sci., 4 (2011), pp.~789--806.

\bibitem{HaScSc-M2AS05}
{\sc M.~Haltmeier, T.~Schuster, and O.~Scherzer}, {\em Filtered backprojection
  for thermoacoustic computed tomography in spherical geometry}, Math. Methods
  Appl. Sci., 28 (2005), pp.~1919--1937.

\bibitem{Hristova-IP09}
{\sc Y.~Hristova}, {\em Time reversal in thermoacoustic tomography - an error
  estimate}, Inverse Problems, 25 (2009).
\newblock 055008.

\bibitem{Isakov-Book02}
{\sc V.~Isakov}, {\em Inverse {Problems} for {Partial} {Differential}
  {Equations}}, Springer-Verlag, New York, second~ed., 2002.

\bibitem{KiSc-SIAM13}
{\sc A.~Kirsch and O.~Scherzer}, {\em Simultaneous reconstructions of
  absorption density and wave speed with photoacoustic measurements}, SIAM J.
  Appl. Math., 72 (2013), pp.~1508--1523.

\bibitem{KuKu-EJAM08}
{\sc P.~Kuchment and L.~Kunyansky}, {\em Mathematics of thermoacoustic
  tomography}, Euro. J. Appl. Math., 19 (2008), pp.~191--224.

\bibitem{KuKu-HMMI10}
\leavevmode\vrule height 2pt depth -1.6pt width 23pt, {\em Mathematics of
  thermoacoustic and photoacoustic tomography}, in Handbook of Mathematical
  Methods in Imaging, O.~Scherzer, ed., Springer-Verlag, 2010, pp.~817--866.

\bibitem{Kunyansky-IP08}
{\sc L.~Kunyansky}, {\em Thermoacoustic tomography with detectors on an open
  curve: an efficient reconstruction algorithm}, Inverse Problems, 24 (2008).
\newblock 055021.

\bibitem{LiWa-PMB09}
{\sc C.~Li and L.~Wang}, {\em Photoacoustic tomography and sensing in
  biomedicine}, Phys. Med. Biol., 54 (2009), pp.~R59--R97.

\bibitem{LiMa-Book72}
{\sc J.-L. Lions and E.~Magenes}, {\em Non-Homogeneous Boundary Value
  Problems}, Springer, Berlin, 1972.

\bibitem{MaRe-CMS14}
{\sc A.~V. Mamonov and K.~Ren}, {\em Quantitative photoacoustic imaging in
  radiative transport regime}, Comm. Math. Sci., 12 (2014), pp.~201--234.

\bibitem{McLean-Book00}
{\sc W.~McLean}, {\em Strongly Elliptic Systems and Boundary Integral
  Equations}, Cambridge University Press, Cambridge, 2000.

\bibitem{NaSc-SIAM14}
{\sc W.~Naetar and O.~Scherzer}, {\em Quantitative photoacoustic tomography
  with piecewise constant material parameters}, SIAM J. Imag. Sci., 7 (2014),
  pp.~1755--1774.

\bibitem{Nguyen-IPI09}
{\sc L.~V. Nguyen}, {\em A family of inversion formulas in thermoacoustic
  tomography}, Inverse Probl. Imaging, 3 (2009), pp.~649--675.

\bibitem{PaSc-IP07}
{\sc S.~K. Patch and O.~Scherzer}, {\em Photo- and thermo- acoustic imaging},
  Inverse Problems, 23 (2007), pp.~S1--S10.

\bibitem{QiStUhZh-SIAM11}
{\sc J.~Qian, P.~Stefanov, G.~Uhlmann, and H.~Zhao}, {\em An efficient
  {Neumann}-series based algorithm for thermoacoustic and photoacoustic
  tomography with variable sound speed}, SIAM J. Imaging Sci., 4 (2011),
  pp.~850--883.

\bibitem{QiSa-SIAM15}
{\sc L.~Qiu and F.~Santosa}, {\em Analysis of the magnetoacoustic tomography
  with magnetic induction}, SIAM J. Imag. Sci., 8 (2015), pp.~2070--2086.

\bibitem{ReGaZh-SIAM13}
{\sc K.~Ren, H.~Gao, and H.~Zhao}, {\em A hybrid reconstruction method for
  quantitative photoacoustic imaging}, SIAM J. Imag. Sci., 6 (2013),
  pp.~32--55.

\bibitem{Scherzer-Book10}
{\sc O.~Scherzer}, {\em Handbook of Mathematical Methods in Imaging},
  Springer-Verlag, 2010.

\bibitem{StUh-IP09}
{\sc P.~Stefanov and G.~Uhlmann}, {\em Thermoacoustic tomography with variable
  sound speed}, Inverse Problems, 25 (2009).
\newblock 075011.

\bibitem{Tittelfitz-IP12}
{\sc J.~Tittelfitz}, {\em Thermoacoustic tomography in elastic media}, Inverse
  Problems, 28 (2012).
\newblock 055004.

\bibitem{TuWe-Book09}
{\sc M.~Tucsnak and G.~Weiss}, {\em Observation and Control for Operator
  Semigroups}, Birkhauser Verlag, Basel, 2009.

\bibitem{Wang-Book09}
{\sc L.~V. Wang}, ed., {\em Photoacoustic Imaging and Spectroscopy}, Taylor \&
  Francis, 2009.

\bibitem{Yamamoto-IP95}
{\sc M.~Yamamoto}, {\em Stability, reconstruction formula and regularization
  for an inverse source hyperbolic problem by a control method}, Inverse
  Probl., 11 (1995), pp.~481--496.

\bibitem{Zuazua-CEAUCT01}
{\sc E.~Zuazua}, {\em Some results and open problems on the controllability of
  linear and semilinear heat equations}, in Carleman Estimates and Applications
  to Uniqueness and Control Theory, F.~Colombini and C.~Zuily, eds.,
  Birkha\"user, Boston, MA, 2001.

\end{thebibliography}

}

\end{document}